\documentclass[]{amsart}
\usepackage[foot]{amsaddr}
\usepackage[a4paper,top=2.5cm,bottom=3cm,inner=3cm,outer=3cm]{geometry}

\usepackage[utf8]{inputenc}

\usepackage{amsmath,amsthm}
\usepackage{amsfonts,amssymb,mathrsfs}

\usepackage{tikz, tikz-cd}
\usetikzlibrary{babel}

\usepackage[bookmarks=false,breaklinks=true]{hyperref}
\usepackage{float}
\usepackage{graphicx}
\usepackage[percent]{overpic}
\usepackage{color}

\usepackage{epigraph}
\definecolor{myurlcolor}{rgb}{0,0,0.7}
\usepackage{hyperref}
\hypersetup{colorlinks,
linkcolor=myurlcolor,
citecolor=myurlcolor,
urlcolor=myurlcolor}
\usepackage[mathscr]{euscript}

\input xy
\xyoption{all}

\newcommand{\maps}{\colon}    

\newcommand{\C}{\mathsf{C}}
\newcommand{\Set}{\mathsf{Set}}
\newcommand{\Vect}{\mathsf{Vect}}

\newcommand{\op}{\mathrm{op}}

\newcommand{\define}[1]{{\bf \boldmath{#1}}}

\theoremstyle{definition}

        \newcommand{\be}{\begin{equation}}
        \newcommand{\ee}{\end{equation}}
        \newcommand{\ba}{\begin{eqnarray}}
        \newcommand{\ea}{\end{eqnarray}}
        \newcommand{\ban}{\begin{eqnarray*}}
        \newcommand{\ean}{\end{eqnarray*}}
        \newcommand{\barr}{\begin{array}}
        \newcommand{\earr}{\end{array}}

\begin{document}
\title{Isbell Duality}
\author[Baez]{John C.\ Baez} 
\address{Department of Mathematics, University of California, Riverside CA, 92521, USA}
\address{Centre for Quantum Technologies, National University of Singapore, 117543, Singapore}
\date{\today}
\maketitle

Mathematicians love dualities.  The dual of a vector space $V$ is the vector space $V^\ast$ of linear maps from $V$ to the ground field.   Any linear map $f \maps V \to W$ between vector spaces gives a linear map going the other way between their duals, $f^\ast \maps W^\ast \to V^\ast$, given by
\[      f^\ast(\ell)(v) = \ell(f(v)) , \qquad \forall v \in V, \ell \in W^\ast.   \]
Composition gets turned around:
\[          (fg)^\ast = g^\ast f^\ast  . \]
Furthermore, there is always a linear map
\[ i \maps  V  \to V^{\ast \ast} \]
given by
\[      i(v)(\ell) = \ell(v)    \qquad \forall v \in V, \ell \in V^*, \]
and when $V$ is finite-dimensional this is an isomorphism.   So, for finite-dimensional vector spaces, duality is like flipping a coin upside down: when you do it twice, you get back where you started---at least up to isomorphism.  
 
Dualities are useful because they let you view the same situation in two  different ways.    Often dualities give an interesting description of the opposite of a familiar category $\C$.  This is the category $\C^{\op}$ with the same objects as $\C$, but where a morphism $f \maps X \to Y$ is defined to be a morphism $f \maps Y \to X$ in $\C$, and the order of composition is reversed.    

Some dualities show a category is equivalent to its own opposite.  For example, duality for vector spaces can be used to show the category of finite-dimensional vector spaces over any field is equivalent to its opposite.   We're secretly using this whenever we take the transpose of a matrix. Similarly, Pontryagin duality says the category of locally compact abelian groups is equivalent to its own opposite.   This duality sends each locally compact group $G$ to a locally compact group called its `Pontryagin dual' $\widehat{G}$, and the Fourier transform of a function on $G$ is a function on $\widehat{G}$.   For example, the Poincar\'e dual of the real line is the real line, while the Pontryagin dual of the circle is $\mathbb{Z}$.

More commonly, however, a category of mathematical objects is not equivalent to its opposite, and a duality relates two different categories.  The opposite of a category of spaces is typically a category of commutative rings or algebras.   For example, Gelfand--Naimark duality says the opposite of the category of compact Hausdorff spaces is the category of commutative $C^\ast$-algebras.  In fact, to make the duality between spaces and commutative rings as nice as possible, Grothendieck \emph{defined} a category of spaces called `affine schemes' to be the opposite of category of commutative rings.


There are also dualities within category theory itself.  The opposite of a category is itself a kind of
a dual, and taking the opposite twice gives you back the category you started with:
\[     (\C^\op)^\op = \C .\]
But there is a subtler and very beautiful duality in category theory called `Isbell duality'.   

First, there is a map from any category $\C$ to the category $[\C^\op,\Set]$, where objects are functors from $\C^{\op}$ to the category of sets and morphisms are natural transformations.   This map takes any object $X \in \C$ to the functor 
\[                \hom(-, X) \maps \C^\op \to \Set \]
which sends any object $A \in \C$ to the set $\hom(A,X)$ of all morphisms from $A$ to $X$.  
This map is called the \define{Yoneda embedding}, and is itself a functor:
\[      \begin{array}{rccl}            y \maps &\C &\to& [\C^\op,\Set]  \\
                                                            &   X & \mapsto & \hom(-,X) .
\end{array}
\]                                                          

The Yoneda embedding is fundamental in category theory.  Philosophically it says that an object can be known by the behavior of the morphisms into it.   It takes time to learn how to use it as a practical tool, but this is nicely explained in modern textbooks \cite{Leinster,Riehl}.  One insight is this: just as we often take a set and form the vector space with that set as basis, it is often useful to treat a category $\C$ as sitting inside the larger category $[\C^\op, \Set]$.   The reason is that $[\C^\op, \Set]$ has `colimits', which are analogous to linear combinations in a vector space.  For example, we can sum $F, G \in [\C^\op,\Set]$ as follows:
\[    (F+G)(X) = F(X) + G(X)   \qquad \forall X \in \C^{\op} \]
where the sum at right is the usual disjoint union of sets.  In fact $[\C^\op, \Set]$ is the free 
category with colimits on the category $\C$.

But this whole story has a dual version!   An object can also be known by the behavior of morphisms \emph{out of} it.  This fact is captured by the \define{co-Yoneda embedding}:
\[      \begin{array}{rccl}            z \maps &\C &\to& [\C,\Set]^{\op}  \\
                                                            &   X & \mapsto & \hom(X,-) .
\end{array}
\]     
The concept dual to colimit is `limit': just as colimits generalize sums, limits generalize products. Unsurprisingly, it turns out that $[\C,\Set]^{\op}$ is the free category with limits on the category $\C$. 

In 1960, Isbell \cite{Isbell} noticed a wonderful link between the Yoneda and co-Yoneda embeddings, which has subsequently been clarified by many authors, as reviewed in \cite{AveryLeinster}.   Any functor $F \maps \C^\op \to \Set$ has an \define{Isbell conjugate} $F^* \maps \C \to \Set$, given by
\[     F^*(X) = \hom(F , y(X)). \]
Similarly, any functor $G \maps \C \to \Set$ has an Isbell conjugate $G^* \maps \C^\op \to \Set$ given by
\[     G^*(X) = \hom(z(X), G) .\]
These two versions of Isbell conjugate give functors going back and forth like this:
\[
\begin{tikzcd}
   \lbrack \C^\op , \Set \rbrack
    \arrow[r, bend left, "\textstyle \ast", pos = 0.5]
    &
    \lbrack \C, \Set \rbrack^\op
    \arrow[l, bend left, "\textstyle \ast", pos = 0.5]
\end{tikzcd}\]
But these two functors are typically not inverses, not even up to natural isomorphism!   Instead, Isbell duality says they are \define{adjoints}, meaning that $\hom(F^*,G)$ and $\hom(F,G^*)$ are naturally isomorphic for all $F \in [\C^{\op},\Set]$ and $G \in [\C, \Set]^{\op}$.  This is analogous to the situation for vector spaces that are not necessarily finite-dimensional: taking the dual defines adjoint functors going back and forth between $\Vect$ and $\Vect^{\op}$, where $\Vect$ is the category of \emph{all} vector spaces over a given field.

Isbell duality sets the stage for a panoply of further developments.  For example, just as the vector spaces where $i \maps V \to V^{\ast \ast}$ is an isomorphism are precisely the finite-dimensional ones, it is very interesting to study functors $F \in [\C^\op,\Set]$ such that the canonical map $i \maps F \to F^{\ast \ast}$ is an isomorphism.  These objects form a subcategory of $[\C^\op,\Set]$ called the \define{reflexive completion} of $\C$ \cite{AveryLeinster,HughesPavlovic}.

So far most applications of Isbell duality involve its generalization to enriched
categories \cite{Willerton}.   For example, this generalization gives a way to find, for any compact metric space $X$, the smallest compact metric space $Y$ containing a copy of $X$ with the property that for any $x \in X$ and $y \in Y$ there is a point $x' \in X$ such that
\[          d(x,y) + d(y,x') = d(x,x') . \]
However, it seems that Isbell duality still remains largely unexploited.    Perhaps one problem is simply that this jewel of mathematics is not yet widely known.


\begin{thebibliography}{99}  

\bibitem{AveryLeinster} T.\ Avery and T.\ Leinster, Isbell conjugacy and the reflexive completion, \textsl{Theory Appl.\ Categ.\ } \textbf{36} (2021), 306--347.  Also available as \href{https://arxiv.org/abs/2102.08290}{arXiv:2102.08290}.

\bibitem{HughesPavlovic} D.\ J.\ D.\ Hughes and D.\ Pavlovic,
The nucleus of an adjunction and the Street monad on monads.  Available as
\href{https://arxiv.org/abs/2004.07353}{arXiv:2004.07353}.

\bibitem{Isbell} J.\ R.\ Isbell, Adequate subcategories, \textsl{Illinois Jour.\ Math.\ } \textbf{4} (1960), 541--552.

\bibitem{Leinster} T.\ Leinster, \textsl{Basic Category Theory}, Cambridge U.\ Press, Cambridge, 2014.    Also available as \href{https://arxiv.org/abs/1612.09375}{arXiv:1612.09375}.


\bibitem{Riehl} E.\ Riehl, \textsl{Category Theory in Context}, Dover, New York, 2016.  Also available as \href{http://www.math.jhu.edu/~eriehl/context.pdf}{http://www.math.jhu.edu/$\sim$eriehl/context.pdf}

\bibitem{Willerton} S.\ Willerton, Tight spans, Isbell completions and semi-tropical modules.  Available as \href{https://arxiv.org/abs/1302.4370}{arXiv:1302.4370}

\end{thebibliography}
\end{document}